\documentclass[11pt,a4paper]{amsart}
\usepackage{amssymb,amsmath}
\usepackage{latexsym}
\usepackage{amsthm,amsfonts,amssymb,mathrsfs}
\usepackage{rotating}
\usepackage[leqno]{amsmath}
\usepackage{xspace}
\usepackage[all]{xy}
\usepackage{longtable}
\usepackage{amstext}
\usepackage{amscd}
\usepackage{enumerate}
\usepackage{setspace}
\usepackage{mathrsfs}

\headsep=1cm \numberwithin{equation}{section}
\newtheorem{theorem}{Theorem}[section]
\newtheorem{lemma}[theorem]{Lemma}

\newtheorem{corollary}[theorem]{Corollary}

\theoremstyle{definition}
\newtheorem{definition}[theorem]{Definition}
\theoremstyle{remark}
\newtheorem{remark}[theorem]{Remark}
\newtheorem{example}[theorem]{Example}
\newtheorem{question}[theorem]{Question}

\newcommand{\Ass}{\operatorname{Ass}}

\newcommand{\grade}{\operatorname{grade}}

\newcommand{\Spec}{\operatorname{Spec}}

\newcommand{\RH}{\operatorname{H}}

\newcommand{\Ht}{\operatorname{ht}}
\newcommand{\id}{\operatorname{id}}

\newcommand{\pd}{\operatorname{pd}}

\newcommand{\E}{\operatorname{E}}

\newcommand{\Ext}{\operatorname{Ext}}
\newcommand{\Supp}{\operatorname{Supp}}

\newcommand{\Tor}{\operatorname{Tor}}
\newcommand{\Hom}{\operatorname{Hom}}

\newcommand{\Ann}{\operatorname{Ann}}

\newcommand{\depth}{\operatorname{depth}}

\newcommand{\Coker}{\operatorname{Coker}}

\newcommand{\lo}{\longrightarrow}
\newcommand{\fm}{\mathfrak{m}}
\newcommand{\fp}{\mathfrak{p}}

\newcommand{\fa}{\mathfrak{a}}

\newcommand{\fn}{\mathfrak{n}}

\newenvironment{prf}[1][Proof]{\begin{proof}[\bf #1]}{\end{proof}}

\begin{document}

\author[K. Abolfath Beigi, K. Divaani-Aazar and M. Tousi]{Kosar Abolfath Beigi, Kamran
Divaani-Aazar and Massoud Tousi, Tehran}

\title[On the invariance of certain types ...]
{On the invariance of certain types of generalized Cohen-Macaulay modules under Foxby equivalence}

\address{K. Abolfath Beigi, Department of Mathematics, Faculty of Mathematical Sciences, Alzahra
University, Tehran, Iran.}
\email{kosarabolfath@gmail.com}

\address{K. Divaani-Aazar, Department of Mathematics, Faculty of Mathematical Sciences, Alzahra
University, Tehran,Iran -and-School of Mathematics, Institute for Research in Fundamental Sciences
(IPM), P.O. Box 19395-5746, Tehran, Iran.}
\email{kdivaani@ipm.ir}

\address{M. Tousi, Department of Mathematics, Faculty of Mathematical Sciences, Shahid Beheshti
University, Tehran, Iran.}
\email{mtousi@ipm.ir}

\begin{abstract} Let $R$ be a local ring and $C$ a semidualizing module of $R$. We investigate the behavior
of certain classes of generalized Cohen-Macaulay $R$-modules under the Foxby equivalence between the Auslander
and Bass classes with respect to $C$. In particular, we show that generalized Cohen-Macaulay $R$-modules are
invariant under this equivalence and if $M$ is a finitely generated $R$-module in the Auslander class with
respect to $C$ such that $C\otimes_RM$ is surjective Buchsbaum, then $M$ is also surjective Buchsbaum.
\end{abstract}

\maketitle

\section{Introduction}

Throughout this paper, $\left(R,\fm,k\right)$ is a commutative Noetherian local ring with nonzero identity.
Let $\mathscr{P}\left(R\right)$ and (resp. $\mathscr{I}\left(R\right)$) denote the full subcategory of finitely
generated $R$-modules of finite projective (resp. injective) dimension.

Let $R$ be Cohen-Macaulay with a dualizing module $\omega_R$. By virtue of \cite[Theorem 2.9]{Sh}, there is an
equivalence of categories

\begin{displaymath}
\xymatrix{\mathscr{P}\left(R\right) \ar@<0.7ex>[rrr]^-{\omega_R\otimes_R-} & {} & {} & \mathscr{I}\left(R\right).
\ar@<0.7ex>[lll]^-{\Hom_R\left(\omega_R,-\right)}}
\end{displaymath}
Let $\mathcal{CM}\left(R\right)$ denote the full subcategory of  Cohen-Macaulay $R$-modules. By \cite[Theorem 3.3(i)]{K},
the above equivalence leads to the following one

\begin{displaymath}
\xymatrix{\mathscr{P}\left(R\right)\cap \mathcal{CM}\left(R\right) \ar@<0.7ex>[rrr]^-{\omega_R\otimes_R-} & {} & {} &
\mathscr{I}\left(R\right)\cap \mathcal{CM}\left(R\right).
\ar@<0.7ex>[lll]^-{\Hom_R\left(\omega_R,-\right)}}
\end{displaymath}

The notion of generalized Cohen-Macaulay modules is one of the most natural extensions of the notion of Cohen-Macaulay
modules. Surjective Buchsbaum modules and Buchsbaum modules are special classes of generalized Cohen-Macaulay modules.
Kawasaki \cite[Theorem 3.3(iii)]{K} showed that if an $R$-module $M\in \mathscr{P}\left(R\right)$ is such that the
$R$-module $\omega_R\otimes_RM$ is surjective Buchsbaum, then $M$ itself is so. Also, he presented an example to show
that the converse doesn't hold; see \cite[Proposition 4.1]{K}.

Let $\left(R,\fm,k\right)$ be an arbitrary local ring and $C$ a semidualizing module of $R$. Let $\mathscr{A}_C\left(R\right)$ and $\mathscr{B}_C\left(R\right)$ denote the Auslander and Bass classes with respect to $C$; respectively. It is known that $\mathscr{P}\left(R\right)\subseteq \mathscr{A}_C\left(R\right)$ and $\mathscr{I}\left(R\right)\subseteq \mathscr{B}_C\left(R\right)$. Also,
these classes are closely related to Gorenstein homological dimensions. Namely, if $R$ is Cohen-Macaulay with a dualizing module $\omega_R$,
then $\mathscr{A}_{\omega_R}\left(R\right)$ (resp. $\mathscr{B}_{\omega_R}\left(R\right)$) is precisely the class of all $R$-modules with
finite Gorenstein projective (resp. Gorenstein injective) dimension; see \cite[Theorems 4.1 and 4.4]{CFH}. By Foxby equivalence, there is an equivalence of categories:

\begin{displaymath}
\xymatrix{\mathscr{A}_C\left(R\right) \ar@<0.7ex>[rrr]^-{C\otimes_R-} &
{} & {} & \mathscr{B}_C\left(R\right).  \ar@<0.7ex>[lll]^-{\Hom_R\left(C,-\right)}}
\end{displaymath}

In this article, we investigate the behavior of the classes of Gorenstein, maximal Cohen-Macaulay, Cohen-Macaulay, surjective Buchsbaum, Buchsbaum and generalized Cohen-Macaulay $R$-modules under the later equivalence. More precisely, if $\mathscr{X}
\left(R\right)$ is any of these classes of modules, we examine the existence of an equivalence
\begin{displaymath}
\xymatrix{\mathscr{A}_C\left(R\right)\cap \mathscr{X}\left(R\right) \ar@<0.7ex>[rrr]^-{C\otimes_R-} &
{} & {} & \mathscr{B}_C\left(R\right)\cap \mathscr{X}\left(R\right).  \ar@<0.7ex>[lll]^-{\Hom_R\left(C,-\right)}}
\end{displaymath}

To ease the reading of the paper, in Section 2, we will recall the needed definitions and bring some elementary results.
In Section 3, we will consider the classes of Gorenstein, maximal Cohen-Macaulay and Cohen-Macaulay $R$-modules; see
Theorem \ref{33}.  In Section 4, which is the core of this paper, we will consider the classes of surjective Buchsbaum,
Buchsbaum and generalized Cohen-Macaulay $R$-modules; see Theorems \ref{43}, \ref{47} and Corollary \ref{49}. Also,
we provide an example of a Cohen-Macaulay local ring $R$ with a dualizing module $\omega_R$ and a Buchsbaum $R$-module
$M\in \mathscr{A}_{\omega_R}(R)$ such that $\omega_R\otimes_RM$ is not Buchsbaum; see Example \ref{48}.

\section{Prerequisites}

We begin with recalling the definitions of weak sequences and filter regular sequences.

\begin{definition} Let $M$ be a finitely generated $R$-module and $\underline{x}=x_1, \dots, x_n$ be a sequence
of elements of $\fm$.
\begin{itemize}
\item[(i)] Following \cite[page 39]{T}, we say $\underline{x}$ is a {\it weak $M$-sequence} if $$\left(\langle x_{1}, x_2,
\ldots, x_{i-1}\rangle M \underset{M}:x_i \right)\subseteq \left(\langle x_{1}, x_2, \ldots, x_{i-1} \rangle M \underset{M}:
\fm \right)$$ for every $i=1,\ldots, n$.
\item[(ii)] Following \cite{CST}, we say $\underline{x}$ is a {\it filter $M$-regular sequence} if $$\left(\langle x_{1}, x_2,
\ldots, x_{i-1}\rangle M \underset{M}:x_i \right)\subseteq \bigcup\limits_{\ell\in \mathbb{N}} \left(\langle x_{1}, x_2, \ldots,
x_{i-1} \rangle M \underset{M}:\fm^{\ell} \right)$$ for every $i=1,\ldots, n$.
\end{itemize}
\end{definition}

By the definition, it is clear that any weak $M$-sequence is a filter $M$-regular sequence.

The following lemma is well-known and easily verified.

\begin{lemma}\label{22} Let $M$ be a finitely generated $R$-module and $\underline{x}=x_1,\dots, x_n\in \fm$. The following are equivalent:
\begin{itemize}
\item[(i)] $\underline{x}$ is a filter $M$-regular sequence.
\item[(ii)] For each $1\leq i\leq r$, the element $x_{i}$ doesn't belong to the union of the elements of $\Ass_R\left(M/\langle x_{1}, x_2,
\ldots, x_{i-1}\rangle M \right)\setminus \{\fm\}.$
\item[(iii)] $x_1/1,\dots, x_n/1$ is a weak $M_\fp$-regular sequence for every $\fp \in \Supp_R(M)\setminus\lbrace \fm\rbrace$.
\end{itemize}
\end{lemma}

In the sequel, we denote the length of an $R$-module $M$ by $\ell_R\left(M\right)$.

\begin{definition} Let $M$ be a finitely generated $R$-module.
\begin{itemize}
\item[(i)] We say that $M$ is {\it Gorenstein} if either $M=0$ or $M$ is non-zero and $\id_RM=\depth_RM$.
\item[(ii)] We say that $M$ is {\it maximal Cohen-Macaulay} if $\depth_RM=\dim R$.
\item[(iii)] We say that $M$ is {\it surjective Buchsbaum} if the natural map $$\phi^i_M :\Ext^i_R\left(k,M\right)\longrightarrow
\RH^i_{\fm}\left(M\right)$$ is surjective for every $i<\dim_RM$.
\item[(iv)] We say that $M$ is {\it Buchsbaum} if every system of parameters $x_1,\dots, x_d$ of $M$ forms a weak $M$-sequence.
\item[(v)] We say that $M$ is {\it generalized Cohen-Macaulay} if $\ell_R\left(\RH^i _\fm\left(M\right)\right)<\infty$ for every
$i<\dim_RM$.
\end{itemize}
\end{definition}

\begin{remark}\label{24}
\begin{itemize}
\item[(i)] Every Cohen-Macaulay $R$-module is surjective Buchsbaum.
\item[(ii)] The classes of modules introduced in the above definition are subset of each other from top to bottom; respectively.
\item[(iii)] If $M$ is a generalized Cohen-Macaulay $R$-module, then every system of parameters of $M$ is a filter $M$-regular
sequence; see \cite[Appendix, Remark 11 and Theorem 14]{SV}.
\item[(iv)] If $R$ is a quotient of a Cohen-Macaulay local ring, then an $R$-module $M$ is generalized Cohen-Macaulay
if and only if every system of parameters of $M$ is a filter $M$-regular sequence; see \cite[Appendix, Proposition 16]{SV}.
\end{itemize}
\end{remark}

\begin{definition} A finitely generated $R$-module $C$ is called {\it semidualizing} if it satisfies the following conditions:
\begin{itemize}
\item[(i)] the homothety map $\chi_C^R:R \lo \Hom_R\left(C,C\right)$ is an isomorphism; and
\item[(ii)] $\Ext^i_R\left(C,C\right)=0$ for all $i>0$.
\end{itemize}
\end{definition}

\begin{definition} Let $C$ be a semidualizing module of $R$.
\begin{itemize}
\item[(i)] The Auslander class $\mathscr{A}_C\left(R\right)$ is the class of all $R$-modules $M$ for which the natural map
$\gamma_M^C:M\lo \Hom_R\left(C,C\otimes_RM\right)$ is an isomorphism and $$\Tor^R_i\left(C,M\right)=0=\Ext_R^i\left(C,C\otimes_RM\right)$$ for all
$i\geq 1$.
\item[(ii)] The Bass class $\mathscr{B}_C\left(R\right)$ is the class of all $R$-modules $M$ for which the evaluation map
$\xi_M^C:C\otimes_R\Hom_R\left(C,M\right)\lo M$ is an isomorphism and $$\Ext^i_R\left(C,M\right)=0=\Tor^R_i\left(C,\Hom_R\left(C,M\right)\right)$$ for all
$i\geq 1$.
\end{itemize}
\end{definition}

\section{Cohen-Macaulay modules}

In Theorem \ref{33}, we show that  Cohen-Macaulay $R$-modules are invariant under Foxby equivalence. To prove it, we need the following two lemmas.

\begin{lemma}\label{31} Let $C$ be a semidualizing module of $R$ and $M$ a finitely generated $R$-module. Then:
\begin{itemize}
\item[(i)] $\Supp_RC=\Spec R$.
\item[(ii)] $\Supp_R\left(C\otimes_RM\right)=\Supp_RM=\Supp_R\left(\Hom_R\left(C,M\right)\right)$. In particular, $$\dim_R\left(C\otimes_RM\right)=\dim_RM=
\dim_R\left(\Hom_R\left(C,M\right)\right).$$
\item[(iii)] $C\otimes_RM\neq 0$ if and only if $M\neq 0$ if and only if $\Hom_R\left(C,M\right)\neq 0$.
\end{itemize}
\end{lemma}

\begin{prf} (i) It is well-known and obvious,  because $$\Spec R=\Supp_R\left(\Hom _R \left(C,C\right)\right)\subseteq \Supp _RC\subseteq
\Spec R.$$

(ii) By (i), one has $$\Supp_R \left(C\otimes_RM\right)=\Supp_RM\cap \Supp_RC=\Supp_RM.$$ On the other hand, \cite[page 267,
Proposition 10]{B} implies that $$\Ass_R\left(\Hom_R\left(C,M\right)\right)=\Supp_RC\cap \Ass_RM=\Ass_RM,$$ and so $\Supp_RM=
\Supp_R\left(\Hom_R\left(C,M\right)\right)$.

(iii) is clear by (ii).
\end{prf}

\begin{lemma}\label{32} Let $C$ be a semidualizing module of $R$ and $\underline{x}=x_1, \dots, x_n\in \fm$. Assume that
$M\in \mathscr{A}_{C}\left(R\right)$ and $N\in \mathscr{B}_{C}\left(R\right)$ are two finitely generated $R$-modules. Then:
\begin{itemize}
\item[(i)] $\Ass_RM=\Ass_R\left(C\otimes_RM\right)$ and $\Ann_RM=\Ann_R\left(C\otimes_RM\right)$.
\item[(ii)]  $\underline{x}$ is a (weak) $M$-regular sequence if and only if it is a (weak) $C\otimes_RM$-regular sequence.
In particular, $\depth_RM=\depth_R\left(C\otimes_RM\right)$.
\item[(iii)] $\Ass_RN=\Ass_R\left(\Hom_R\left(C,N\right)\right)$ and $\Ann_RN=\Ann_R\left(\Hom_R\left(C,N\right)\right)$.
\item[(iv)]  $\underline{x}$ is a (weak) $N$-regular sequence if and only if it is a (weak) $\Hom_R\left(C,N\right)$-regular
sequence. In particular, $\depth_RN=\depth_R\left(\Hom_R\left(C,N\right)\right)$.
\end{itemize}
\end{lemma}

\begin{prf} (i) By \cite[p. 267, Proposition 10]{B}, one has
$$\begin{array}{ll}
\Ass_RM&=\Ass_R\left(\Hom_R\left(C,C\otimes_RM\right)\right)\\
&=\Supp_RC\cap \Ass_R\left(C\otimes_RM\right)\\
&=\Ass_R\left(C\otimes_RM\right).
\end{array}
$$
By considering the $R$-isomorphism $M\cong \Hom_R\left(C,C\otimes_RM\right)$, one can easily deduce that $\Ann_RM=
\Ann_R\left(C\otimes_RM\right)$.

(ii)  By Lemma \ref{31}(iii), it follows that $M/\underline{x}M \neq 0$ if and only if $\left(C\otimes_RM\right)/\underline{x}\left(C\otimes_RM\right)\neq 0$. So, it is enough
to show that $\underline{x}$ is a weak $M$-regular sequence if and only if it is a weak
$C\otimes_RM$-regular sequence.  We prove this by induction on $n$. The case $n=1$ follows
by (i). Note that for a finitely generated $R$-module $L$, the union of members of $\Ass_RL$
is the set of all elements of $R$ which are zero-divisor on $L$.

Next, assume that $n>1$ and the claim holds for $n-1$. By the case $n=1$, without loss of generality, we can and do assume
that $x_1$ is regular on both $M$ and $C\otimes_RM$. Set $\overline{M}:=M/x_1M$. Applying \cite[Proposition 3.1.7(a)]{Sa} on
the exact sequence $$0 \longrightarrow M \overset{x_1}\longrightarrow  M \longrightarrow\overline{M} \longrightarrow 0,$$ yields
that $\overline{M}\in \mathscr{A}_C\left(R\right)$. By induction hypothesis, $x_2, \dots, x_n$ is a weak $\overline{M}$-regular
sequence if and only if $x_2, \dots, x_n$ is a weak $C\otimes_R\overline{M}$-regular sequence. Since  $$C\otimes_R\overline{M}\cong \left(C\otimes_RM\right)/x_1\left(C\otimes_RM\right),$$ the proof is complete.

As $N\in \mathscr{B}_{C}\left(R\right)$, one has $N\cong C\otimes_R\Hom_R\left(C,N\right)$ and $\Hom_R\left(C,N\right)\in \mathscr{A}_{C}\left(R\right)$. Thus (iii) and (iv) are immediate by (i) and (ii); respectively.
\end{prf}

Part (i) of the next result extends \cite[Theorem 3.3(i)]{K}.

\begin{theorem}\label{33} Let $C$ be a semidualizing module of $R$ and $M$ a finitely generated $R$-module in $\mathscr{A}_C\left(R\right)$.
Then:
\begin{itemize}
\item[(i)] $M$ is Cohen-Macaulay if and only if $C\otimes_RM$ is Cohen-Macaulay.
\item[(ii)] $M$ is maximal Cohen-Macaulay if and only if $C\otimes_RM$ is maximal Cohen-Macaulay.
\item[(iii)] If $M$ is Gorenstein, then $C\otimes_RM$ is also Gorenstein.
\end{itemize}
\end{theorem}

\begin{prf} (i) By Lemma \ref{31}(ii) and Lemma \ref{32}(ii), one has $\dim_RM=\dim_R\left(C\otimes_RM\right)$ and $\depth_RM=
\depth_R\left(C\otimes_RM\right)$. So, (i) is immediate.

(ii) It is immediate, because by Lemma \ref{32}(ii) we have $\depth_RM=\depth_R\left(C\otimes_RM\right)$.

(iii) Assume that $M$ is Gorenstein and set $N:=C\otimes_RM$. If $M=0$, then the claim is obvious. So, we assume that $M\neq 0$. Then $N\neq 0$
and $N\in \mathscr{B}_C\left(R\right)$. As $N\neq 0$, $\depth_RN<\infty$. Now, as $\Ext_R^i(C,N)=0$ for all $i\geq 1$ and $\Hom_R\left(C,N\right)
\cong M$ has finite injective dimension, \cite[Theorem 8.1]{C} yields that $C\cong R$. Thus $C\otimes_RM\cong M$, and so $C\otimes_RM$ is Gorenstein.
\end{prf}

The following example indicates that the converse of Theorem \ref{33}(iii) doesn't hold.

\begin{example}\label{34} Let $\left(R,\fm,k\right)$ be a Cohen-Macaulay local ring which is not Gorenstein with a dualizing module $\omega_R$.
Then $\omega_R\otimes_RR$ is Gorenstein, but $R$ is not Gorenstein.
\end{example}

\section{Generalized Cohen-Macaulay modules}

In this section, we first prove that generalized Cohen-Macaulay $R$-modules are invariant under Foxby equivalence. To this end, we need the following
two lemmas.

\begin{lemma}\label{41} Let $C$ be a semidualizing module of $R$, $M$ an $n$-dimensional finitely generated $R$-module and $\underline{x}
=x_1, \dots, x_n\in \fm$. Then $\underline{x}$ is a system of parameters of $M$ if and only if it is a system of parameters of $C\otimes_RM$.
\end{lemma}

\begin{prf} First note that by Lemma \ref{31}(ii), we have $\dim_R\left(C\otimes_RM\right)=\dim_RM$. Applying Lemma \ref{31}(ii) again yields that
$$\begin{array}{ll}
\dim_R\left(\left(C\otimes_RM\right)/\langle x_{1}, x_2, \ldots, x_{n}\rangle \left(C\otimes_RM\right)\right)&=\dim_R
\left(C\otimes_R\left(M/\langle x_{1}, x_2, \ldots, x_{n}\rangle M\right)\right)\\
&=\dim_R\left(M/ \langle x_{1}, x_2, \ldots, x_{n} \rangle M\right).
\end{array}
$$
This completes the argument, because for an $n$-dimensional finitely generated $R$-module $N$, it is known that $\underline{x}$
is a system of parameters of $N$ if and only if $\dim_R\left(N/\langle \underline{x} \rangle N\right)=0$.
\end{prf}

\begin{lemma}\label{42} Let $C$ be a semidualizing module of $R$. Let $M$ be a finitely generated $R$-module in $\mathscr{A}_C\left(R\right)$ and
$\underline{x}=x_1, \dots, x_n\in \fm$. Then $\underline{x}$ is a filter $M$-regular sequence if and only if $\underline{x}$ is a
filter $C\otimes_RM$-regular sequence.
\end{lemma}

\begin{prf} Let $\fp\in \Spec R$. Then by \cite[Propositions 2.2.3 and 3.5.3]{Sa}, respectfully, $C_{\fp}$ is a semidualizing
module of the local ring $R_\fp$ and $M_\fp\in \mathscr{A}_{C_\fp}\left(R_\fp\right)$. Thus Lemma \ref{32}(ii) yields that $x_1/1,\dots, x_n/1$
is a weak $M_\fp$-regular sequence if and only if $x_1/1,\dots, x_n/1$ is a weak $C_{\fp}\otimes_{R_{\fp}}M_{\fp}$-regular sequence. Now,
Lemma \ref{22} completes the proof. Note that $C_{\fp}\otimes_{R_{\fp}}M_{\fp}\cong {\left(C\otimes_RM\right)}_{\fp}$ and by Lemma \ref{31}(ii), $\Supp_RM=\Supp_R\left(C\otimes_RM\right)$.
\end{prf}

The next result generalizes \cite[Theorem 3.3(ii)]{K}.

\begin{theorem}\label{43} Let $C$ be a semidualizing module of $R$ and $M$ a finitely generated $R$-module in $\mathscr{A}_C\left(R\right)$.
Then $M$ is generalized Cohen-Macaulay if and only if $C\otimes_RM$ is generalized Cohen-Macaulay.
\end{theorem}

\begin{prf} By \cite[Propositions 2.2.1 and 3.4.7]{Sa}, $\widehat{C}$ is a semidualizing module of $\widehat{R}$ and $\widehat{M}\in
\mathscr{A}_{\widehat{C}}\left(\widehat{R}\right)$. On the other hand, it is routine to check that a finitely generated $R$-module $N$
is generalized Cohen-Macaulay if and only if the $\widehat{R}$-module $\widehat{N}$ is so. Also, there is a natural $\widehat{R}$-isomorphism
$\widehat{C}\otimes_{\widehat{R}}\widehat{M}\cong \widehat{C\otimes_RM}$. Thus, we may and do assume that $R$ is complete. Now, the assertion follows by Lemmas \ref{41} and \ref{42} and Remark \ref{24}(iv). Note that by the Cohen Structure Theorem, every complete local ring is a homomorphic image of a regular local ring.
\end{prf}

Let $\fa$ be an ideal of $R$ and $M$ and $N$ two $R$-modules. The $i$th generalized local cohomology module of $M$ and $N$ with respect to $\fa$ is defined by $\RH^i_{\fa}(M,N):=\underset{n}{\varinjlim}\Ext^{i}_{R}(M/\fa^{n}M,N)$; see \cite{H}.

\begin{lemma}\label{44} Let $\fa$ be an ideal of $R$ and $C$ a semidualizing module of $R$. Then for every $R$-module $N$ in $\mathscr{B}_C\left(R\right)$,
there is a natural $R$-isomorphism $\RH_{\fa}^i\left(C,N\right)\cong \RH_{\fa}^i\left(\Hom_R\left(C,N\right)\right)$ for all $i\geq 0$.
\end{lemma}

\begin{prf} For every finitely generated $R$-module $L$ and every injective $R$-module $E$, we claim that $\RH_{\fa}^i\left(\Hom_R\left(L,E\right)\right)=0$ for all $i\geq 1$. To this end, let $F_{\bullet}$ be a
free resolution of $L$ consisting of finitely generated free $R$-modules. As the $R$-module $E$ is injective,
it follows that $\Hom_R\left(F_{\bullet},E\right)$ is an injective resolution of the $R$-module $\Hom_R\left(L,E\right)$.
Since $\Gamma_{\fa}(E)$ is an injective $R$-module, for each natural integer $i$, one has $$\RH_{\fa}^i\left(\Hom_R\left(L,E\right)\right)=\RH^i\left(\Gamma_{\fa}\left(\Hom_R\left(F_{\bullet},E\right)\right)\right)=
\RH^i\left(\Hom_R\left(F_{\bullet},\Gamma_{\fa}\left(E\right)\right)\right)=0.$$

Let $I^{\bullet}$ be an injective resolution of $N$. Since $N$ and all modules in the complex $I^{\bullet}$ belong to $\mathscr{B}_C\left(R\right)$, it follows that the functor	$\Hom_R(C,-)$ leaves the exact sequence $$0\lo N\lo I^0\lo
\cdots \lo I^j\lo \cdots $$ exact. Hence, $\Hom_R\left(C,I^{\bullet}\right)$ is a $\Gamma_{\fa}$-acyclic resolution of
the $R$-module $\Hom_R\left(C,N\right)$. This implies the first isomorphism in the following
display:
$$\begin{array}{ll}
\RH_{\fa}^i\left(\Hom_R\left(C,N\right)\right)&\cong \RH^i\left(\Gamma_{\fa}\left(\Hom_R\left(C,I^{\bullet}\right)\right)\right)\\
&\cong \RH^i\left(\underset{n}{\varinjlim}\Hom_R\left(R/\fa^n,\Hom_R\left(C,I^{\bullet}\right)\right)\right)\\
&\cong \underset{n}{\varinjlim}\RH^i\left(\Hom_R\left(R/\fa^n,\Hom_R\left(C,I^{\bullet}\right)\right)\right)\\
&\cong \underset{n}{\varinjlim}\RH^i\left(\Hom_R\left(C/\fa^n C,I^{\bullet}\right)\right)\\
&\cong \underset{n}{\varinjlim}\Ext^{i}_{R}\left(C/\fa^{n}C,N\right)\\
&=\RH_{\fa}^i\left(C,N\right).
\end{array}
$$
\end{prf}

\begin{corollary}\label{45} Let $C$ be a semidualizing module of $R$ and $M$ a finitely generated $R$-module in $\mathscr{A}_C\left(R\right)$.
Then there is a spectral sequence $$\Ext^p_{\fm} \left(C,\RH^{q}_{\fm} \left(C\otimes_RM\right)\right) \underset{p}\Longrightarrow
\RH^{p+q}_{\fm}\left(M\right).$$
\end{corollary}

\begin{prf} Let $L$ and $N$ be two finitely generated $R$-modules. Consider the covariant functors $G(-):=\Gamma_{\fm}(-)$ and
$F(-):=\Hom_R(L,-)$ on the category of $R$-modules and $R$-homomorphisms. Then, for every injective $R$-module $I$, the $R$-module
$G(I)$ is injective and $\RH^{i}_{\fm}\left(L,N\right)\cong R^i(FG)(N)$ for all $i\in \mathbb{N}_0$. Hence, applying \cite[Theorem 10.47]{R} on the functors $G(-)$ and $F(-)$ yields a spectral sequence $$\Ext^p_R \left(L,\RH^q_{\fm}\left(N\right)\right)\underset{p}\Longrightarrow \RH^{p+q}_{\fm}\left(L,N\right).$$ Letting $L:=C$ and $N:=C\otimes_RM$, by Lemma \ref{44}, one deduces that $$\RH^{p+q}_{\fm}\left(C,C\otimes_RM\right)\cong \RH^{p+q}_{\fm}\left(\Hom_R\left(C,C\otimes_RM\right)\right)\cong
\RH^{p+q}_{\fm}\left(M\right).$$
\end{prf}

Recall that for each non-negative integer $n$,  the $n$th Betti (resp. Bass) number of an $R$-module $M$ is defined as $\beta_n^{R}\left(M\right):=
\text{Vdim}_k\left(\Tor^n_R\left(k,M\right)\right)$ (resp.  $\mu^n_{R}\left(M\right):=\text{Vdim}_k\left(\Ext^n_R\left(k,M\right)\right)$).
If there is no ambiguity about the underlying ring, we then show these invariants with $\beta_n\left(M\right)$ and $\mu^n\left(M\right)$; respectively.

\begin{lemma}\label{46} Let $M$ be a generalized Cohen-Macaulay $R$-module with $d:=\dim_RM>0$. Then we have $$\mu^n\left(M\right)\leqslant \sum^{n}_{j=0}\beta_j\left(k\right)\ell_R\left(\RH^{n-j}_{\fm}\left(M\right)\right)$$ for all $n\leq d$. Furthermore, the following are equivalent:
\begin{itemize}
\item[(i)] $M$ is a surjective Buchsbaum $R$-module.
\item[(ii)] The equality in above holds for all $n<d$.
\end{itemize}
\end{lemma}

\begin{prf} See \cite[Corollary 1.14]{M}, \cite[Theorem 1.2]{Y} and \cite[Lemma 2.2]{K}.
\end{prf}

The next result is a far reach generalization of \cite[Theorem 3.3(iii)]{K}.

\begin{theorem}\label{47} Let $C$ be a semidualizing module of $R$ and $M$ a finitely generated $R$-module in $\mathscr{A}_C\left(R\right)$.
Assume that $C\otimes_RM$ is surjective Buchsbaum. Then $M$ is also surjective Buchsbaum.
\end{theorem}

\begin{prf} Remark \ref{24}(ii) and Theorem \ref{43} imply that $M$ is generalized Cohen-Macaulay. Clearly, we may assume that $d:=\dim_RM>0$.
So by Lemma \ref{46}, we have $$\mu^n\left(M\right)\leqslant\sum \limits^n_{j=0} \beta_j\left(k\right) \ell_R\left(\RH^{n-j}_{\fm}\left(M\right)\right)
\  \ \left(*\right)$$ and $$\mu^{n}\left(C\otimes_RM\right)=\sum\limits^n _{j=0}\beta_j \left(k\right)\ell_R\left(\RH^{n-j}_{\fm}\left(C\otimes_RM\right)
\right)   \  \  \left(**\right)$$ for all $n<d$. Note that by Lemma \ref{31}(ii), $\dim_RM=\dim_R\left(C\otimes_RM\right).$ Let $0\leqslant n < d$ be an
integer. Consider the spectral sequence $$\Ext^p_R \left(C,\RH^q _{\fm}\left(C\otimes_RM\right)\right)\underset{p}\Longrightarrow \RH^{p+q} _{\fm}
\left(M\right);$$ see Corollary \ref{45}.
There is a filtration $$0=\RH^{-1}\subseteq \RH^0\subseteq\dots\subseteq \RH^{n-1}\subseteq \RH^n=\RH^n _{\fm}\left(M\right)$$
such that $\RH^p/\RH^{p-1}\cong \E^{p,n-p}_\infty$ for all $p=0,\ldots, n$. Hence, $$\ell_R\left(\RH^n_{\fm}\left(M\right)\right)=\sum^n_{p=0}
\ell_R\left(\E^{p,n-p}_\infty \right)\leqslant \sum^n_{p=0} \ell_R\left(\E^{p,n-p}_2\right).$$ Note that $\E^{p,n-p}_\infty$ is a subquotient
of $\E^{p,n-p}_2$
for all $p=0,\ldots, n.$ Let $F_{\bullet}$ be a minimal free resolution of $C$. Then $\E^{p,n-p}_2=\Ext^p_R\left(C,\RH^{n-p}_{\fm}
\left(C\otimes_RM\right)\right)$ is
a subquotient of $\Hom_R\left(F_p,\RH^{n-p}_{\fm}\left(C\otimes_RM\right)\right)$, and so $$\ell_R \left(\RH^n_{\fm}\left(M\right)\right)
\leqslant \sum^n_{p=0}\beta_p\left(C\right)
\ell_R \left(\RH^{n-p}_{\fm}\left(C\otimes_RM\right)\right).  \  \ \left(\dag\right)$$

As $M\cong \Hom_R\left(C,C\otimes_RM\right)$ and $\Ext_R^i\left(C,C\otimes_RM\right)=0$ for all $i\geq 1$,  we get $M\simeq R\Hom_R
\left(C,C\otimes_RM\right)$,
and so \cite[13.19]{F} yields that $$\mu^n \left(M\right)=\sum^n_{p=0}\beta_p\left(C\right)
\mu^{n-p}\left(C\otimes_RM\right).$$ Thus,
$$\begin{array}{ll}
\mu^n\left(M\right)&\overset{\left(**\right)}=\sum\limits_{p=0}^{n}\beta_p\left(C\right)\left(\sum\limits_{j=0}^{n-p}\beta_j\left(k\right) \ell_R\left(\RH^{n-p-j}_{\fm}\left(C\otimes_RM\right)\right)\right)\\
&=\sum\limits^{n}_{j=0}\beta_j\left(k\right)\left(\sum\limits^{n-j}_{p=0} \beta_p\left(C\right) \ell_R\left(\RH^{n-j-p}_{\fm}\left(C\otimes_RM
\right)\right)\right)\\
&\overset{\left(\dag\right)}\geqslant \sum\limits_{j=0}^{n} \beta_j\left(k\right)\ell_R\left(\RH^{n-j}_{\fm}\left(M\right)\right)\\
&\overset{\left(*\right)}\geqslant \mu^n \left(M\right).
\end{array}
$$
So, $$\mu^n\left(M\right)=\sum^n_{j=0} \beta_j \left(k\right) \ell_R \left(\RH^{n-j}_{\fm}\left(M\right)\right)$$ for all $n<d$. Therefore,
$M$ is surjective Buchsbaum by
Lemma \ref{46}.
\end{prf}

Next, we provide an example of a Cohen-Macaulay local ring $R$ with a dualizing module $\omega_R$ and a Buchsbaum $R$-module $M$ of finite
projective dimension such that the $R$-module $\omega_R\otimes_RM$ is not Buchsbaum.

\begin{example}\label{48} Let $\left(T,\fn,k\right)$ be a 4-dimensional regular local ring with $k$ infinite. Let $\fa\subset {\fn}^2$ be
an ideal of $T$ with 3 generators such that $R:=T/\fa$ is a 2-dimensional Cohen-Macaulay local ring which is not Gorenstein. (For an explicit
realization,
let $k$ be an infinite field, $T:=k[[X,Y,Z,W]]$ and $R:=k[[X,Y,Z,W]]/\langle X^4-Y^3,X^5-Z^3,Y^5-Z^4 \rangle$.)

The Auslander-Buchsbaum formula yields that $\pd_TR=2$. Let $\mu\left(\fa\right)$ denote the minimum number of generators of $\fa$. Then
$$2=\dim T-\dim R=\text{ht} \fa\leq \mu\left(\fa\right)\leq 3.$$ We claim that $\mu\left(\fa\right)=3.$ In the contrary, suppose that $\fa$
can be generated by two elements
$z_1$ and
$z_2$. Then by \cite[Corollarry 1.6.19]{BH}, $z_1, z_2$ is a $T$-regular sequence, and so $R$ is Gorenstein. From this contradiction, we conclude
that $\mu\left(\fa\right)=3.$

As $\mu\left(\fa\right)=3$, the minimal free resolution of the $T$-module $R$ has the form $$0\lo T^n\lo T^3\lo T\lo 0.$$ We show that $n=2$.
To this end, let $K$ denote the quotient field of the domain $T$. As $\fa \neq 0$, we immediately see that $K\otimes_T\fa\cong K$. Hence, by
applying the exact functor $K\otimes_T-$ to the exact sequence $$0\lo T^n\lo T^3\lo \fa\lo 0,$$ it turns out that $n=2$.

Let $\omega_R$ denote the dualizing module of $R$. Then by \cite[Corollarry 3.3.9]{BH}, the minimal free resolution of the $T$-module $\omega_R$
has the form $$0\lo T\lo T^3\lo T^2\lo 0.$$ Applying the right exact functor $R\otimes_T -$ to the exact sequence $T^3\lo T^2\lo \omega_R\lo 0$, implies the exact sequence $R^3\lo R^2\lo \omega_R\lo 0$. Thus $R^3\lo R^2\lo 0$ is the beginning of the minimal free resolution of the $R$-module
$\omega_R$, and so $\beta^R_0\left(\omega_R\right)=\beta^T_0\left(\omega_R\right)=2$ and $\beta^R_1\left(\omega_R\right)=\beta^T_1\left(\omega_R\right)=3$.

By \cite[Proposition 4.1]{K} and its proof, there exists a 2-dimensional surjective Buchsbaum $R$-module $M$ of finite projective dimension such that the
$R$-module $\omega_R\otimes_RM$ is not surjective Buchsbaum and

\begin{equation}
\ell_T\left(\RH^i_{\fn}\left(\omega_R\otimes_RM\right)\right)=
\left\{\begin{array}{ll}
\beta^T_1\left(R\right) & i=1;\\
\beta^T_0\left(R\right) & i=0;
\end{array} \right.
\end{equation}
and
\begin{equation}
\ell_T\left(\RH^i_{\fn}\left(M\right)\right)=
\left\{\begin{array}{ll}
\beta^T_1\left(\omega_R\right) & i=1;\\
\beta^T_0 \left(\omega_R\right)& i=0.
\end{array} \right.
\end{equation}

We claim that the $R$-module $\omega_R\otimes_RM$ is not even Buchsbaum. Suppose the contrary holds. Then by \cite[Chapter I, Lemma 1.6]{SV}, $M$ and
$\omega_R\otimes_RM$ are also Buchsbaum as $T$-modules. As $T$ is regular, \cite[Chapter I, Corollary 2.16]{SV} yields that $M$ and $\omega_R\otimes_RM$
are surjective Buchsbaum $T$-modules. In particular, by Lemma \ref{46}, one has the following equalities:

\begin{equation}
\begin{array}{ll}
\mu ^1_T \left(\omega_R\otimes_RM\right)&=\beta^T_0\left(k\right)\ell_T\left(\RH^1_{\fn}\left(\omega_R\otimes_RM\right)\right)+\beta^T_1\left(k\right)
\ell_T\left(\RH^0_{\fn}\left(\omega_R\otimes_RM\right)\right)\\
&=\beta^T_0\left(k\right)\beta^T_1\left(R\right)+\beta^T_1\left(k\right)\beta^T_0\left(R\right)\\
&=3+\beta^T_1\left(k\right).
\end{array}
\end{equation}
and
\begin{equation}
\begin{array}{ll}
\mu ^1_T \left(M\right)&=\beta^T_0\left(k\right)\ell_T\left(\RH^1_{\fn}\left(M\right)\right)+\beta^T_1\left(k\right)\ell_T\left(\RH^0_{\fn}\left(M\right)\right)\\
&=\beta^T_0\left(k\right)\beta^T_1\left(\omega_R\right)+\beta^T_1\left(k\right)\beta^T_0\left(\omega_R\right)\\
&=3+2\beta^T_1\left(k\right).
\end{array}
\end{equation}

Since $\pd_RM$ is finite, it follows that $M\in \mathscr{A}_{\omega_R}\left(R\right)$. In particular, $\Tor^R_i(\omega_R,M)=0$ for all $i\geq 0$, and so one has
the quism $\omega_R\otimes_RM\simeq \omega_R\otimes_R^LM$. By \cite[Exercise  3.3.26]{BH}, we have $\mu^1_T\left(\omega_R\otimes_RM\right)=\beta^T_3\left(\omega_R
\otimes_RM\right)$ and $\mu^1_T\left(M\right)=\beta^T_3\left(M\right)$. So,
\begin{equation}
\begin{array}{ll}
\mu^1_T\left(\omega_R\otimes_RM\right)&=\text{Vdim}_k\left(\text{H}_3\left(\left(\omega_R\otimes_RM\right)\otimes_T^Lk\right)\right)\\
& =\text{Vdim}_k\left(\text{H}_3\left(\left(\omega_R\otimes_R^LM\right)\otimes_T^Lk\right)\right)\\
& =\text{Vdim}_k\left(\text{H}_3\left(\omega_R\otimes_R^L\left(M\otimes_T^Lk \right)\right)\right)\\
&=\sum \limits_{i+j=3}\text{Vdim}_k\left(\text{H}_i\left(\omega_R\otimes_R^Lk\right)\right)\times\text{Vdim}_k\left(\text{H}_j
\left(M\otimes_T^Lk\right)\right)\\
& \geq \text{Vdim}_k\left(\text{H}_0\left(\omega_R\otimes_R^Lk\right)\right)\times \text{Vdim}_k\left(\text{H}_3\left(M\otimes_T^Lk\right)\right)\\
& =\beta^R_0\left(\omega_R\right)\beta^T_3\left(M\right)\\
& =2\mu^1_T\left(M\right).
\end{array}
\end{equation}

The fourth equality holds by \cite[10.7]{F}. Now, (4.3) and (4.4) yield the desired contradiction. Hence, the $R$-module $\omega_R\otimes_RM$ is not Buchsbaum.
\end{example}

\begin{corollary}\label{49} Let $C$ be a semidualizing module of $R$ and $M$ a finitely generated $R$-module in $\mathscr{B}_C\left(R\right)$.
\begin{itemize}
\item[(i)] $M$ is generalized Cohen-Macaulay if and only if $\Hom_R\left(C,M\right)$ is generalized Cohen-Macaulay.
\item[(ii)] If $M$ is surjective Buchsbaum, then $\Hom_R\left(C,M\right)$ is also  surjective Buchsbaum.
\item[(iii)] $M$ is Cohen-Macaulay if and only if $\Hom_R\left(C,M\right)$ is Cohen-Macaulay.
\item[(iv)] $M$ is maximal Cohen-Macaulay if and only if $\Hom_R\left(C,M\right)$ is maximal Cohen-Macaulay.
\item[(v)] If the $R$-module $\Hom_R\left(C,M\right)$ is Gorenstein, then $M$ is also Gorenstein.
\end{itemize}
\end{corollary}

\begin{prf} In view of the natural isomorphism $M\cong C\otimes_R\Hom_R\left(C,M\right)$ and the fact $\Hom_R(C,M)\in \mathscr{A}_C\left(R\right)$,
these claims follow immediately by Theorems \ref{33}, \ref{43} and \ref{47}.
\end{prf}

Let $C$ be a semidualizing module of $R$ and $n$ be a non-negative integer.  Let  $\mathscr{P}_C\left(R\right)$ and
$\mathscr{I}_C\left(R\right)$ denote, respectively, the classes of $C$-projective and $C$-injective finitely generated
$R$-modules. Set $\widehat{\mathscr{P}_C}\left(R\right)_{\leq n}$, $\widehat{\mathscr{P}}\left(R\right)_{\leq n}$,
$\widehat{\mathscr{I}_C}\left(R\right)_{\leq n}$, and $\widehat{\mathscr{I}}\left(R\right)_{\leq n}$ to be the classes
of finitely generated $R$-modules of $C$-projective, projective, $C$-injective and injective dimension of at most $n$;
respectively. Then by \cite[Theorem 2.12]{TW}, there are the following Foxby equivalences of categories:

\begin{displaymath}
\xymatrix@C=20ex{\mathscr{P}(R) \ar@{^(->}[d] \ar@<0.8ex>[r]_-{}^-{C\otimes_R-} & \mathscr{P}_C(R)\ar@{^(->}[d]\ar@<0.8ex>[l]^-{\Hom_R(C,-)} \\
\widehat{\mathscr{P}}(R)_{\leq n} \ar@{^(->}[d] \ar@<0.8ex>[r]_-{}^-{C\otimes_R-} & \widehat{\mathscr{P}}_C(R)_{\leq n} \ar@{^(->}[d]\ar@<0.8ex>[l]^-{\Hom_R(C,-)} \\
\mathscr{A}_C(R) \ar@<0.8ex>[r]_-{}^-{C \otimes_R-} & \mathscr{B}_C(R) \ar@<0.8ex>[l]^-{\Hom_R(C,-)} \\
\widehat{\mathscr{I}}_C(R)_{\leq n} \ar@{^(->}[u] \ar@<0.8ex>[r]_-{}^-{C \otimes_R-} & \widehat{\mathscr{I}}(R)_{\leq n}. \ar@<0.8ex>[l]^-{\Hom_R(C,-)} \ar@{^(->}[u] \\
\mathscr{I}_C(R) \ar@{^(->}[u] \ar@<0.8ex>[r]_-{}^-{C \otimes_R-} & \mathscr{I}(R). \ar@<0.8ex>[l]^-{\Hom_R(C,-)} \ar@{^(->}[u]}
\end{displaymath}

In view of Theorems \ref{33}, \ref{43} and \ref{47} and Corollary \ref{49}, we deduce the following corollary.

\begin{corollary}\label{410} Let $C$ be a semidualizing module of $R$ and $n$ a non-negative integer. Let $\left(\mathscr{Y}\left(R\right),\mathscr{Z}\left(R\right)\right)$
be any of the pairs $\left(\mathscr{P}\left(R\right),\mathscr{P}_C\left(R\right)\right)$, $\left(\widehat{\mathscr{P}}\left(R\right)_{\leq n},\widehat{\mathscr{P}}_C\left(R\right)_{\leq n}\right)$,
$\left(\mathscr{A}_C\left(R\right),\mathscr{B}_C\left(R\right)\right)$, $\left(\widehat{\mathscr{I}}_C\left(R\right)_{\leq n},\widehat{\mathscr{I}}\left(R\right)_{\leq n}\right)$; or $\left(\mathscr{I}_C\left(R\right),
\mathscr{I}\left(R\right)\right)$. Then there is the following display of equivalences and functors:
\begin{displaymath}
\xymatrix@C=20ex{
\mathscr{Y}(R)\cap \mathcal{G}(R) \ar@{^(->}[d] \ar@<0.8ex>[r]_-{}^-{C\otimes_R-} & \mathscr{Z}(R)\cap \mathcal{G}(R)\ar@{^(->}[d] \\
\mathscr{Y}(R)\cap \mathcal{MCM}(R) \ar@{^(->}[d] \ar@<0.8ex>[r]_-{}^-{C\otimes_R-} & \mathscr{Z}(R)\cap \mathcal{MCM}(R) \ar@{^(->}[d] \ar@<0.8ex>[l]^-{\Hom_R(C,-)} \\
\mathscr{Y}(R)\cap \mathcal{CM}(R) \ar@{^(->}[d]\ar@<0.8ex>[r]_-{}^-{C \otimes_R-} & \mathscr{Z}(R)\cap \mathcal{CM}(R) \ar@{^(->}[d] \ar@<0.8ex>[l]^-{\Hom_R(C,-)}\\
\mathscr{Y}(R)\cap \mathcal{SB}(R) \ar@{^(->}[d] & \mathscr{Z}(R)\cap \mathcal{SB}(R) \ar@{^(->}[d]\ar@<0.8ex>[l]^-{\Hom_R(C,-)}\\
\mathscr{Y}(R)\cap \mathcal{GCM}(R) \ar@<0.8ex>[r]_-{}^-{C \otimes_R-} & \mathscr{Z}(R)\cap \mathcal{GCM}(R). \ar@<0.8ex>[l]^-{\Hom_R(C,-)}}
\end{displaymath}
\end{corollary}

Here, $\mathcal{G}\left(R\right)$ and $\mathcal{MCM}\left(R\right)$ denote the classes of Gorenstein and maximal Cohen-Macaulay
$R$-modules; respectively and $\mathcal{SB}$ and $\mathcal{GCM}$ stand for the full subcategory of surjective Buchsbaum and generalized Cohen-Macaulay $R$-modules; respectively.

We end the paper by proposing the following natural question.

\begin{question}\label{411} Let $C$ be a semidualizing module of $R$ and $M$ a finitely generated $R$-module in $\mathscr{A}_C\left(R\right)$
such that $C\otimes_RM$ is Buchsbaum. Then, is $M$ Buchsbaum?
\end{question}


\end{document}